\documentclass{article}

\usepackage{amsmath}
\usepackage{amssymb}
\usepackage[latin1]{inputenc} 
\usepackage{graphicx}
\usepackage{epsfig}
\usepackage{color}
\usepackage{float}
\usepackage[all]{xy}
\usepackage{gastex}
\newcommand{\n}{{\noindent}}
\renewcommand{\b}{{\bigskip}}
\newcommand{\q}{{\quad}}                 

\renewcommand{\r}{\right}
     
\newtheorem{theorem}{Theorem}
\newtheorem{lemma}{Lemma}

\newtheorem{corollary}{Corollary}

\begin{document}

\title{Frises}
\author{ Ibrahim Assem \and Christophe Reutenauer \and David Smith }
\maketitle

\begin{abstract}
Each acyclic graph, and more generally, each acyclic orientation of the graph associated to a Cartan matrix, allows to define a so-called frise; this is a collection of sequences over $\mathbf N$, one for each vertex of the graph. We prove that if these sequences satisfy a linear recurrence, then the Cartan matrix is of Dynkin type (if the sequences are bounded) or of Euclidean type (if the sequences are unbounded). We prove the converse in all cases, except for the exceptional Euclidean Cartan matrices; we show even that the sequences are rational over the semiring $\mathbf N$. We generalize these results by considering frises with variables; as a byproduct we obtain, for the cases case $A_m$ and $\tilde A_m$, explicit formulas for the cluster variables, over the semiring of Laurent polynomials over $\mathbf N$ generated by the initial variables (which explains simultaneously positivity and the Laurent phenomenon). The general tool are the so-called $SL_2$-tilings of the plane; these are fillings of the whole discrete plane by elements of a ring, in such a way that each $2\times 2$ connected submatrix is of determinant 1.
\end{abstract}
\section{Introduction}

To each acyclic quiver (that is, to each acyclic directed graph) and more generally, to each Cartan matrix), we associate its so-called {\em frise}. This construction is motivated by the theory of cluster algebras of Fomin-Zelevinsky, and was shown to two of the authors by Philippe Caldero some years ago. For each such graph, this construction gives sequences of positive rational numbers, one for each vertex of the graph. Actually, the Laurent phenomenon of Fomin-Zelevinsky implies that these sequences are integer-valued. 

The present work was motivated by a simple example, shown in Section 2: the quiver with two vertices and two edges from one to the other. Then the numbers which appears in the sequences of the frise are the Fibonacci numbers of even rank. This lead us to the question of determining those graphs for which the sequences are {\em  rational} (that is, satisfy some linear recursion, as do the Fibonacci numbers). We conjecture that these graphs are exactly the Dynkin and Euclidean (also called affine) graphs, with some orientation.

We show in Theorem 1 that this conjecture holds in one direction: if a Cartan matrix gives rise to a frise with rational sequences, then this Cartan matrix is of Dynkin or Euclidean type. For this, we apply criteria of Vinberg and Berman-Moody-Wonenburger, which characterize these graphs combinatorially, through the so-called additive and subadditive functions. In order to obtain such functions, we take, very roughly speaking, the logarithm of the sequences of the frise (there are some technicalities, due to the fact that rational sequences do not grow very smoothly in general).

We prove the opposite direction of the conjecture for all but a finite number of cases (Theorem 2): for each Dynkin and for each non-exceptional Euclidean graph, the sequences of the frise are rational. Note that the Dynkin case is immediately settled as a consequence of the finite-type classification of cluster algebras by Fomin and Zelevinsky \cite{FZ2}. So only remains to prove the case of Euclidean graphs.

The rationality for the Euclidean graphs is proved by using a mathematical object, which seems interesting in itself: the {\em $SL_2$-tilings} of the plane. This is a filling of the whole discrete plane by elements of a commutative ring, in such a way that each 2 by 2 connected submatrix is of determinant 1. Note that an analogue mathematical object has been already considered in the literature: the so-called "frieze patterns" in \cite{Co}, \cite{Coco1}, \cite{Coco2}. In our settings, the latter would be called partial $SL_2$-tilings, since they do not cover the whole plane.

We give a way to construct $SL_2$-tilings, where the whole plane is filled with natural integers. As a corollary, we obtain Theorem 2 for the case $\tilde A_m$. The case $\tilde D_m$ is more involved, since we have to study special tilings containing perfect squares. The other Euclidean non-exceptional cases may be reduced to the two previous ones.

At the end of the article, we generalize all this by considering frises and $SL_2$-tilings with variables. This is directly motivated by the theory of cluster algebras, and particularly by the search for formulas expressing cluster variables. 

Recall that cluster algebras were introduced by Fomin and Zelevinsky \cite{FZ1}, \cite{FZ2} in order to explain the connection between the canonical basis of a quantised enveloping algebra and total positivity for algebraic groups.  Since then, they turned out to have important ramifications in several fields of mathematics.  Roughly speaking, a cluster algebra is an integral domain with a possibly infinite family of distinguished generators (called cluster variables) grouped into (overlapping) clusters of the same finite cardinality and computed recursively from an initial cluster.  By construction, every cluster variable can be uniquely expressed as a rational function of the elements of any given cluster.  The Laurent phenomenon, established in \cite{FZ1}, asserts that these rational functions are in fact Laurent polynomials with integral coefficients.  Moreover, it was conjectured by Fomin and Zelevinsky that these coefficients are positive.  This latter conjecture is known as the positivity conjecture. 

Our motivation in the last part of the present article was to develop computational tools to derive direct formulas for the cluster variables without going through the recursive process.  We give explicit and simple formulas, involving only matrix products, for all cluster variables (or all but finitely many cluster variables) for cluster algebras 
without coefficients 
of Dynkin type $A_m$ (or Euclidean type $\tilde A_m$, respectively), explaining simultaneously  the Laurent phenomenon and the positivity for initial acyclic clusters.

It has to be noted that the positivity conjecture, and direct formulas for computing the cluster variables, have been obtained in special cases only (see, for instance, \cite{FZ2, CZ, MS, S, S1, CR, CK, MSW, MuPro}). 
In particular, in \cite{CK, CR}, the positivity conjecture was
established for cluster algebras having no coefficients and an acyclic initial
cluster, by using Euler-Poincar\'{e} characteristics
of appropriate Grassmannians of quiver representations. Note also that summation formulas (where the sum is over perfect matchings of certain graphs) for cluster variables appear in \cite{MSW}, over the semiring of Laurent polynomials over $\mathbf N$; these formulas are valid for a wide class of cluster algebras. 

Our paper provides, through a novel approach leading to much simpler formulas, a new elementary proof of the positivity conjecture for all Dynkin diagrams of type $A$ and all Euclidean diagrams of type $A$. Our approach uses, instead of summation formulas, only products of 2 by 2 matrices over the Laurent polynomial semiring over $\mathbf N$.

For type $A_n$, we use partial $SL_2$-tilings, that are equivalent to the frieze patterns of Coxeter and Conway \cite{Co, Coco1, Coco2}. Their importance for cluster algebras of this type has already been noted by Caldero and Chapoton \cite{Cacha}, see also \cite{BaMaTho}. Our approach provides matrix product formulas to compute them, which seem to be new (however, see \cite{Pro}, where a similar matrix approach is used to compute frieze patterns giving the Markoff numbers).

For details about cluster algebras and cluster categories, we refer the reader to \cite{FZ1, FZ2, BMRRT}, and for information concerning linear recurrences and rational series, we refer to \cite{EPSW, BR}.

{\bf Acknowledgments}: the authors thank Anissa Amroun for computer calculations which helped to determine many frises, and Vestislav Apostolov, François Bergeron, Yann Bugeaud, Christophe Hohlweg, Ralf Schiffler and Lauren Williams for discussions and mail exchanges which helped to clarify several points leading to the present article. Special thanks to Gregg Musiker, who indicated to us the notion of frieze patterns in the literature.

\section{An introductory example}

Given a quiver $Q$ (in other words, a directed graph with possibly multiple edges), which we assume to be acyclic, let $V$ be its set of vertices. Define for each $v$ in $V$ a sequence $v(n)$ by the initial condition $v(0)=1$ and the recursion
$$
v(n+1)=\frac{1}{v(n)} (1+\prod_{v\rightarrow w}w(n) \prod_{w\rightarrow v}w(n+1)).
$$
The fact that these equations define uniquely the sequences $v(n)$ follows from the acyclicity of the graph. 

\begin{figure}[htbp]
\begin{center}

 \scalebox{1}{\input{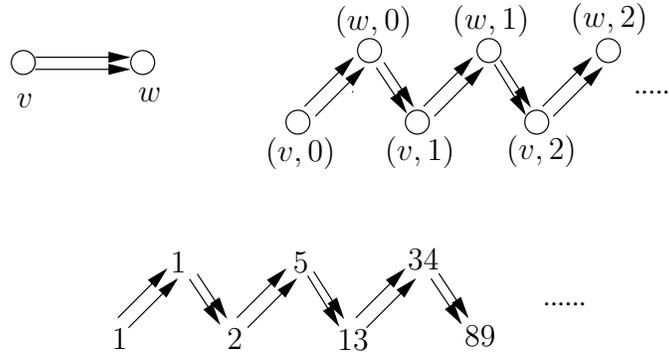}}

\caption{Kronecker quiver and frise}
\label{Kronecker}
\end{center}
\end{figure}

The previous recursion formula  may be represented by defining the {\em frise } associated to the quiver: it is an infinite graph with set of vertices $V\times \mathbb N$ and edges $(v,n)\rightarrow (w,n)$ if $v\rightarrow w$ is in $Q$, and edges $(v,n)\rightarrow (w,n+1)$ if $v\leftarrow w$  in $Q$. Then the sequence $v(n)$ labels the vertex $(v,n)$ of the frise, $l(v,n)=v(n)$ say, and the recursion reads
$$
l(v,n+1)=\frac{1}{l(v,n)}(1+\prod_{(w,i)\rightarrow (v,n+1)} l(w,i)),
$$
with the initial conditions $l(v,0)=1$. Note that only $i=n$ or $i=n+1$ may occur in the product.
 As an example take the {\ Kronecker quiver}, with two vertices and two edges from one to the other. The frise is represented in Figure \ref{Kronecker}, together with the labels. 
 
It turns out that the numbers $1,1,2,5,13,34,...$  are the Fibonacci numbers of even rank  $F_{2n}$, if one defines $F_0=F_1=1$ and $F_{n+2}=F_{n+1}+F_n$. This may be proved for instance by using the identity
$$
\left | \begin{array}{ll} F_{2n+4} & F_{2n+2}\\ F_{2n+2}& F_{2n} \end{array} \right | = 1,
$$
which is a consequence of the fact that the Fibonacci numbers of even rank satisfy the recursion $F_{2n+4}=3F_{2n+2}-F_{2n}$, as is well-known.


Actually, as mentioned before, we shall also deal with more general frises: the initial values $v(0)$ will be variables. In the example of the Kronecker quiver, we may take $u_0=a$, $u_1=b$ and the recursion $u_{n+2}=\frac{1+u_{n+1}^2}{u_{n}}$, which shortcuts the frise: $u_0,u_2,u_4,\dots$ label the vertices $(v,0),(v,1),(v,2),\dots$ and $u_1,u_3, u_5\dots$ the vertices $(w,0),(w,1),(w,2)\dots$. Then one has also a linear recursion, generalizing the recursion for Fibonacci numbers of even rank: 
$$u_{n+2}=\frac{a^2+b^2+1}{ab}u_{n+1}-u_n.$$ 
Moreover, $u_n=\frac{1}{a^{n-2}b^{n-1}}(1,b)M^{n-2}\left ( \begin{array}{l} 1  \\ b\end{array} \right )$, where 
$$
M=\left ( \begin{array}{lc}a^2+1 & b \\ b&b^2\end{array} \right ).
$$
A summation formula for $u_n$ has already been given by Caldero and Zelevinsky \cite{CZ} Th.4.1.

For quivers of type $A_m$ and $\tilde A_m$, we shall obtain these kinds of formulas, which explain simultaneously the Laurent phenomenon (the denominator is a monomial) of Fomin and Zelevinsky, and the positivity of the formulas. 

\section{Frises associated to Cartan matrices}
Recall that a {\em Cartan matrix} $C=(C_{ij})_{1\leq i,j\leq d}$ is defined by the following properties:

\n (i) $C_{ij}\in \mathbf Z$;

\n (ii) $C_{ii}=2$;

\n (iii) $C_{ij}\leq 0$, if $i\neq j$;

\n (iv) $C_{ij} \ne 0 \Leftrightarrow C_{ji}\ne 0$.

The {\em  simple graph associated to}  $C$ has set of vertices $\{1,\dots,d\}$ and an undirected edge $\{i,j\}$ if $i\neq j$ and $C_{ij}\neq 0$. The Cartan matrix is completely described by its {\em diagram}, which is the previous graph, with {\em valuations} on it; the couple $(|C_{ij}|,|C_{ji}|)$ is represented on the edge $\{i,j\}$ as follows:

\begin{figure}[htbp]
\begin{center}

 \scalebox{1}{\input{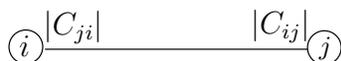}}

\caption{Cartan diagram}
\label{Cartandiagram}
\end{center}
\end{figure}

If $|C_{ij}|=1$, it is omitted. In this manner, Cartan matrices are equivalent to diagrams.

Without loss of generality, we consider only {\em connected} Cartan matrices; this means that the underlying graph is connected. Consider some fixed acyclic orientation of this graph; if the edge $\{i,j\}$ of this graph is oriented from $i$ to $j$, we write $i\rightarrow j$. For each $j=1,...,d$, we define a sequence $a(j,n)$, $n\in \mathbf N$, by the formula, for all $j=1,...,d$:

\begin{equation} \label{frises_induction}
\, a(j,n)a(j,n+1)=1+( \prod_{j\rightarrow i} a(i,n)^{|C_{ij}|})(\prod_{i\rightarrow j} a(i,n+1)^{|C_{ij}|} )
\end{equation}
and the initial conditions $a(j,0)=1$.
 The data of these $d$ sequences is called the {\em frise} associated to the Cartan matrix and the given acyclic orientation of its graph. This generalizes the construction of Section 2. The acyclicity of the orientation ensures that these sequences are well-defined. Moreover, they have clearly coefficients in $\mathbf Q_+^*$. Now, it is a consequence of the Laurent phenomenon of Fomin and Zelevinsky that the coefficients are actually positive integers; see \cite{FZ}; see also \cite{K}. 
 

Recall that a sequence $(a_n)_{n\in \mathbf N}$ of complex numbers {\em satisfies a linear recurrence} if for some $k\geq 1$, some $\alpha_1,...,\alpha_k$ in $\mathbf C$, one has: for all $n$ in $\mathbf N$, $a_{n+k}=\alpha_1 a_{n+k-1}+...+\alpha_k a_n$. Equivalently, the series $\sum_{n\in \mathbf N} a_nx^n$ is {\em rational}, that is, it is the quotient of two polynomials in $\mathbf C[x]$; we say also that the sequence $(a_n)$ is rational. We say that a frise is {\em rational} if the sequences $a(i,n)$ are all rational for $i=1,...,d$.

\b \n \textbf{Main conjecture}. \label{conjecture}
{\em A frise associated to a Cartan matrix with some acyclic orientation is rational if and only if the Cartan matrix is of Dynkin or Euclidean type.}

\n
These Cartan matrices are recalled in Section 9.

One direction of the conjecture is completely solved by the following result.

\begin{theorem} \label{theorem1}
Let $C$ be a connected Cartan matrix with some acyclic orientation. Suppose that the associated frise is rational. If the sequences are all bounded, then $C$ is of Dynkin type. If they are not all bounded, then $C$ is of Euclidean type.
\end{theorem}

\n 
Note that since we assume that the Cartan matrix is connected, the sequences are all simultaneously bounded or unbounded. We prove Theorem 1 in Section 4.

We prove the opposite direction of the conjecture in all cases, except for the exceptional Euclidean cases.
Actually, we prove more than rationality of the sequences. For this, recall that a series $S=\sum_{n\in \mathbf N} a_n x^n \in \mathbf N[[x]]$ is called {\em $\mathbf N$-rational} if it satisfies one of the two equivalent conditions (this equivalence is a particular case of the Kleene-Schützenberger theorem, see \cite{BR}):

\n (i) $S$ belongs to the smallest subsemiring of $\mathbf N[[x]]$ closed under the operation $T\rightarrow T^*=\sum_{n\in \mathbf N} T^n$ (which is defined if $T$ has zero constant term);

\n (ii) for some matrices $\lambda \in \mathbf N^{1\times d}$, $M\in \mathbf N^{d\times d}$, $\gamma \in \mathbf N^{d\times 1}$, one has: for all $n$ in $\mathbf N$, $a_n=\lambda M^n \gamma$. 

We then say that the sequence $(a_n)$ is $\mathbf N$-rational.

A result which we need later in this article is that $\mathbf N$-rational sequences are closed under Hadamard product; this is the coefficientwise product $((a_n),(b_n))\mapsto(a_nb_n)$. This follows easily from (ii) by tensoring the matrices. Note that if $K$ is any commutative semiring, we may define $K$-rational series in exactly the same way, by replacing above $\mathbf N$ by $K$. See \cite{BR}.

\begin{theorem} \label{theorem2}
If a Cartan matrix is of Dynkin or Euclidean type, but not one of the exceptional Euclidean types, then the sequences of each associated frise are $\mathbf N$-rational, and in particular, rational.
\end{theorem}

Theorem 2 will be proved in Section 7. Note that in the Dynkin case, this is an immediate consequence of the finiteness of the set of cluster variables, proved in \cite{FZ2}. Indeed, the construction of the frise is a particular case of the mutations of Fomin and Zelevinsky, by performing mutations only on sources. See for example \cite{K}.


In order to prove Theorem 2, we introduce objects which we call {\em $SL_2$-tiling of the plane}. This is a mapping $t:\mathbf Z^2\mapsto \mathbf N$ such that for any $i,j$ in $\mathbf Z$, $$\left | \begin{array}{ll} t(i,j) & t(i,j+1)\\ t(i+1,j) & t(i+1,j+1) \end{array} \right | = 1.$$
An example is given below; the 1's are boldfaced, since they will play a special role in the sequel; we have not represented the numbers above them. The perfect squares in the tiling will be explained further in the article: see Lemma 4.
$$
\begin{array}{lllllllllllllll}
&&&&&&&&&&\bf 1&\bf 1&\bf 1&\bf 1& \\
&&&&&&&&&\bf 1&\bf 1&2&3&4& \\
&&&&&&&&&\bf 1&2&5&8&11& \\
&&&&&&&&&\bf 1&3&8&13&18& \\
&&&&...&&...&&\bf 1&\bf 1&2^2&11&18&25& ...\\
&&&&&&&&\bf 1&2&9&5^2&41&57& \\
&&&&&&&&\bf1&3&14&39&8^2&89& \\
& & &\bf 1&\bf 1&\bf 1&\bf 1&\bf 1&\bf 1&4&19&53&87&11^2& \\
\bf 1&\bf 1&\bf 1&\bf 1&2&3&4&5&6&5^2&119&332&545&758& \\
\bf 1&2&3&4&9&14&19&24&29&121&24^2&1607&2368&3669&... \\
&&&&&...&&...&&&&...&&& \\
\end{array}
$$
We call {\em frontier} a bi-infinite sequence 
\begin{equation}\label{frontier}
\ldots x_{-3}x_{-2}x_{-1}x_0x_1x_2x_3 \ldots
\end{equation}
with $x_i\in \{x,y\}$. It is called {\em admissible} if there are arbitrarily large and arbitrarily small $i$'s such that $x_i=x$, and arbitrarily large and arbitrarily small $j$'s such that $x_j=y$; in other words, none of the two sequences $(x_n)_{n\geq0}$ and $(x_n)_{n\leq0}$ is ultimately constant. Each frontier may be embedded into the plane: the $x$ (resp. $y$) determine the horizontal (resp. vertical) edges of a bi-infinite discrete path: $x$ (resp. $y$) corresponds to a segment of the form $[(a,b),(a+1,b)]$ (resp $[(a,b),(a,b+1)]$). 
The vertices of the path (that is, the endpoints of the previous segments) get the label 1.



An example is the frontier $$\ldots yxxxyxxxxxyyyxyyyxyxxx\ldots$$
 corresponding to the 1's in the above $SL_2$-tiling.
  
We prove below that an admissible frontier, embedded into the plane, may be extended to a unique $SL_2$-tiling. For this, we need the following notation. Let 
$$M(x)=\left(\begin{array}{cc}1&1\\0&1\end{array}\right) \quad \mbox{and}  \quad M(y)=\left(\begin{array}{cc}1&0\\1&1\end{array}\right).$$

We extend $M$ into a homomorphism from the free monoid generated by $x$ and $y$ into the group $SL_2(\mathbf Z)$.
 
Given an admissible frontier, embedded in the plane as explained previously, let $(u,v)\in \mathbf Z^2$. Then we obtain a finite word, which is a factor of the frontier, by projecting the point $(u,v)$ horizontally and vertically onto the frontier. We call this word the {\em word} of $(u,v)$. It is illustrated in the figure below, where the word of the point $M$ is $yyyxxyx$:
$$
 \begin{array}{ccccccccccccccccc} 
&&&&&&&&&&&&. \\
&&&&&&&&&&&.    \\
&&&&&&&&&&.       \\
&&&&&&&\bf 1&\bf 1&\bf 1 \\
&&&&&\bf 1&\bf 1&\bf 1&|\\
&&&&&\bf 1&&&|\\
&&&&&\bf 1&&&|\\
&&&&\bf 1&\bf 1&-&-&M\\
&&&. \\
&&. \\
&. \\
\end{array}
$$
Note that such a word always begins by $y$ and ends by $x$. We define the word of a point only for points below the frontier; for points above, the situation is symmetric and we omit it.

\begin{theorem}
Given an admissible frontier, there exists a unique $SL_2$-tiling of the plane $t$ extending the embedding of the frontier into the plane. It is defined, for any point $(u,v)$ below the frontier, with associated word $x_1x_2...x_{n+1}$, where $n\geq 1$, $x_i\in\{x,y\}$, by the formula
\begin{equation}\label{tilingformula}
t(u,v)=(1,1)M(x_2)\cdots M(x_n)\left(\begin{array}{c}1\\1\end{array}\right).
\end{equation}
\end{theorem}
Theorem 3 will be proved in Section 5. 
In the formula, note that the first and last letter of the word are omitted. An instance of the formula, for the tiling above, is 
$$
14=(1,1)\left(\begin{array}{cc}1&0\\1&1\end{array}\right)  \left(\begin{array}{cc}1&1\\0&1\end{array}\right) \left(\begin{array}{cc}1&0\\1&1\end{array}\right)^3 \left(\begin{array}{c}1\\1\end{array}\right),
$$
since the word corresponding to 14 in the figure is $y^2xy^3x$.

Theorem 3 will be further generalized, by considering tilings with variables. Indeed, in Section 8, we consider frontiers with variables, instead of 1's as in Theorem 3. Then we show that one obtains an $SL_2$-tiling whose values lie in the semiring of Laurent polynomials over $\mathbf N$ in these variables.



\section{Proof of Theorem 1}

Given two sequences of positive real numbers $(a_k)$ and $(b_k)$, we shall write $a_k\approx b_k$ to express the fact that for some positive constant $C$, one has $\lim_{k\rightarrow \infty}a_k/b_k=C$. 

\begin{lemma}\label{asymptotic}
Let $a(j,n)$, for $j=1,...,d$, be $d$ unbounded sequences of positive integers, each satisfying a linear recurrence. There exist an integer $p\geq 1$, real numbers $\lambda(j,l)\geq 1$ and integers $d(j,l)\geq 0$, for $j=1,...,d$ and $l=0,...,p$, and a strictly increasing sequence $(n_k)_{k\in \mathbf N}$ of nonnegative integers, such that:

\n 
(i) for every  $j=1,...,d$ and every $l=0,...,p, a(j,pn_k+l)\approx\lambda(j,l)^{n_k} n_k^{d(j,l)}$;

\n
(ii) for every  $j=1,...,d$, there exists $l=0,...,p$ such that $\lambda(j,l)>1$ or $d(j,l)\geq 1$;

\n
(iii) for every $j=1,...,d$, $\lambda (j,0)=\lambda(j,p)$ and $d(j,0)=d(j,p)$.
\end{lemma}
\b 
This lemma may be well-known to the specialist of linear recurrent sequences. Since we could not find a precise reference, we give a proof below.

\b
\n\textit{Proof.}

\n Step 1. Recall that each sequence $(a_n)_{n\in \mathbf N}$ satisfying a linear recurrence has a unique expression, called the {\em exponential polynomial}, of the form
\begin{equation} \label{exponential_polynomial}
a_n=\sum_{i=1}^{k} P_i(n)\lambda_i^n,
\end{equation}
for $n$ large enough, where $P_i(n)$ is a nonzero polynomial in $n$ and the $\lambda_i$ are distinct nonzero complex numbers; see e.g. \cite{BR} or \cite{EPSW}. We call the $\lambda_i$ the {\em eigenvalues}
of the sequence $a_n$. The {\em degree} of $\lambda_i$ is $deg(P_i)$. Note that if the $a_n$ are all positive integers, then at least one of its eigenvalues has modulus $|\lambda_i|\geq 1$. The {\em principal part} of the above exponential polynomial is
\begin{equation} \label{principal_part}
n^D\sum_{j}\alpha_j\lambda_j^n,
\end{equation}
where the sum is restricted to those $j$ with $|\lambda_j|$ maximum, $D=deg(P_j)$ maximum for these $j$, and $\alpha_j$ is the coefficient of $n^D$ in $P_j$. We call these $\lambda_j$ the {\em dominating eigenvalues}, {\em dominating modulus} their modulus and $D$ the {\em maximum degree}. Note that if $\lambda_i$ is not a dominating eigenvalue, then either $\lambda_i$ has modulus strictly smaller than the maximum modulus, or its modulus is the maximum modulus, but its degree is strictly smaller than $D$.
For further use, note that if Eq.(\ref{principal_part}) is the principal part of $(a_n)$, then the principal part of  the sequence $a_{n+H}$ is
\begin{equation} \label{principal_part+H}
n^D\sum_{j}\alpha_j\lambda_j^H\lambda_j^n.
\end{equation}

Note also that for any $p\geq 1$ and $l\geq 0$, the eigenvalues of the sequence $(a_{np+l})_{n\in\mathbf N}$ are $p$-th powers of eigenvalues of $(a_n)$: it suffices, to see it, to replace $n$ by $np+l$ in Eq. (\ref{exponential_polynomial}).

\n  Step 2. Consider the subgroup $G$ of $\mathbf C^*$ generated by the eigenvalues of the $d$ sequences $a(j,n)$. It is a finitely generated abelian group and therefore, by the fundamental theorem of finitely generated abelian groups, there exists $p\geq 1$ such that the subgroup $G_1$ generated by the $p$-th powers of any set of elements of $G$ is a free abelian group.

Consider the sequence $(a(j,pn+l))_{n\in\mathbf N}$, with $j=1,...,d$ and $l=0,...,p$. The eigenvalues of the sequence $(a(j,pn+l))_{n\in\mathbf N}$ are $p$-th powers of the eigenvalues of $a_n$ (this follows easily from Eq.(\ref{exponential_polynomial})). Denote by $Z_{j,l}$ the set of its dominating eigenvalues and 
$$Z=\bigcup _{1\geq j\geq d, 0\geq l\geq p} Z_{j,l}.$$
Then $Z\subset G_1$. In particular, no quotient $z/z'$ with $z,z'\in Z$, is a nontrivial root of unity, since $G_1$ is a torsion-free group. Note also that $Z_{j,l}$ is nonempty, since the sequences $a(j,pn+l)$ are positive.

\n  Step 3. Suppose that we have proved the lemma for the $d$ sequences $a'(j,n)=a(j,n+ph)$ for some nonnegative integer $h$. Denote by $\lambda(j,l), d(j,l), (n_k)$ the corresponding numbers and sequences. Then we have
$a'(j,pn_k+l)\approx \lambda(j,l)^{n_k} n_k^{d(j,l)}$. Thus $a(j,p(n_k+h)+l)\approx \lambda(j,l)^{n_k} n_k^{d(j,l)}$. This shows that the lemma then holds for the $d$ sequences $a(j,n)$: we simply replace $(n_k)_{k\in \mathbf N}$ by $(n_k+h)_{k\in \mathbf N}$.
We choose $h$ below.

\n  Step 4. The principal part of the sequence $a(j,pn+l)$ is of the form
\begin{equation}
n^{d(j,l)}\sum_{z\in Z_{j,l}} \alpha_z z^n.   \label{principal_part_a(j,pn+l)}
\end{equation}
Consider the sequence $(\sum_{z\in Z_{j,l}} \alpha_z z^n)_{n\in\mathbf N}$. Since no quotient of distinct elements of $Z_{j,l}$ is a root of unity, we see by the theorem of Skolem-Mahler-Lech (see \cite{BR} Th.4.1 or \cite{EPSW} Th.2.1 ) that the previous sequence has only finitely many zeros. Hence, for some $h$, there is no zero for $n\ge h$. We may choose the same $h$ for each $j=1,...,d$ and $l=0,...,p.$ 

By Step 3 and  Eq.(\ref{principal_part+H}) with $H=ph$, we may therefore assume that the principal part of $a(j,pn+l)$ is Eq. (\ref{principal_part_a(j,pn+l)}), with $\sum_{z\in Z_{j,l}} \alpha_z \neq 0$.

\n  Step 5. Choose some $z(j,l)$ in $Z_{j,l}$. The complex numbers $z/z(j,l)$, for $z\in Z_{j,l}$, $j=1,...,d$ and $l=0,...,p$ have modulus 1 and generate a subroup of $G_1$, which by Step 2 is a finitely generated free abelian subgroup of $\mathbf C^*$.

Observe that if a finite set $E$ of complex numbers of modulus 1 generates a free abelian subgroup of $\mathbf C^*$, then there is a strictly increasing sequence $(n_k)_{k\in \mathbf N}$ of nonnegative integers such that: for every $e\in E, \lim_{k\rightarrow \infty}e^{n_k}=1$. This follows from Kronecker's simultaneous approximation theorem applied to a basis of the previous free abelian group (write the basis as $\exp(2i\pi x)$, for a finite set of $\mathbf Q$-linear independent real numbers $x$, see \cite{HW} Th. 442).

Thus we may assume the existence of $n_k$ with $z^{n_k}\sim_{k\rightarrow\infty} z(j,l)^{n_k}$ for any $j=1,...,d$ and $l=0,...,p$.

\n  Step 6. Going back to the principal part Eq.( \ref{principal_part_a(j,pn+l)}), we see that
\begin{equation*}
n_k^{d(j,l)}  \sum_{z\in Z_{j,l}} \alpha_z z^{n_k} 
\sim_{k\rightarrow \infty} n_k^{d(j,l)} z(j,l)^{n_k} \sum_{z\in Z_{j,l}} \alpha_z
\end{equation*}
since the last sum is nonzero.

Now, since the eigenvalues of $a(j,pn+l)$ which are not in $Z(j,l)$  have modulus strictly smaller than that of $z(j,l)$, or have the same modulus but smaller degree, and since $a(j,n)>0$, we obtain that
\begin{equation*}
a(j,pn_k+l)\approx\lambda(j,l)^{n_k} n_k^{d(j,l)},
\end{equation*}
with $\lambda(j,l)=| z(j,l)|$. 

\n  Step 7. In order to prove (ii), we use the fact that $a(j,n)$ is unbounded; hence, for any $j=1,...,d$, there exists $l=0,...,p-1$ such that $a(pn+l)$ is unbounded. Then at least one of its eigenvalues is of modulus $>1$, and $| z(j,l)|>1$, or otherwise, they have all modulus $\leq 1$ and some $d(j,l)$ must be $>1$.

For (iii), note that $a(j,p(n+1))=a(j,pn+p)$, hence the sequences $a(j,pn)$ and $a(j,pn+p)$ have the same eigenvalues and maximum degrees. Thus  
$ |z(j,0)|=|z(j,p)|$ and $d(j,0)=d(j,p)$.

\begin{flushright}
$\square$
\end{flushright} 

Before proving Theorem 1, we must recall some facts about additive and subadditive functions
of diagrams.
Let $C=(C_{ij})_{1\leq i,j\leq d}$ be a Cartan matrix. An {\em additive} (resp. {\em subadditive})
{\em function for} $C$ is a function $f: \{1,...,d\} \rightarrow \mathbf R_+^{ \ast}$ such that for any $j=1,...,d$, one has $2f(j)=\sum _{i\neq j} f(i)|C_{ij}|$ (resp. $2f(j)\geq \sum _{i\neq j} f(i)|C_{ij}|$). Note that, by the properties of a Cartan matrix, this may equivalently be rewritten as $\sum _{i}f(i)C_{ij} =0$ (resp. $\sum _{i}f(i)C_{ij} \geq 0$). 

The results we need is the following. 

\noindent{\textbf{Theorem}}\label{additive}
{\em A Cartan matrix $C$ is of Euclidean type if and only if there exists an additive function for $C$; it is of Dynkin type if and only if there exists a subadditive function for $C$ which is not additive.}

The second part of this theorem is due to Vinberg \cite{V} and the first to Berman, Moody and Wonenburger \cite{BMW}.
Both results were proved by Happel, Preiser and Ringel \cite{HPR}, under the assumption that the function takes integer values, although this assumption was unnecessary. We need the generalization involving real valued functions, which thus holds by the proof of \cite{HPR}, Theorem p.286.

The idea of the proof of Theorem 1 is as follows: we show first that if the sequences of the frise are rational and bounded, then there exists a subadditive function which is not additive for the diagram. Hence by the theorem above, the diagram is of Dynkin type. The subadditive function is obtained by multiplying $p$ consecutive values of each sequence, $p$ being a common period, and then taking their logarithm. The recurrence relations Eq.(\ref{frises_induction}) imply that this logarithm is a subadditive function which is not additive.

In the case where the sequences are unbounded, the proof is similar by replacing each sequence by its principal part. However, the proof is more technical, since the growth of a rational sequence is in general not of exponential type. Lemma \ref{asymptotic} allows to bypass this difficulty.

\noindent \textbf {Proof of Theorem 1}

\n {\bf Case 1.} The sequences $a(j,n)$, $j=1,...,d$, are all bounded. Since they are integer-valued, they take only finitely many values. Since they satisfy linear recursions, they are ultimately periodic. Let $p$ be a common period and let $n_0$ be such that each sequence is purely periodic for $n\geq n_0$.

Let $b(j)=\prod_{n_0\leq n <n_0+p} a(j,n)$. Note that $b(j)>1$. Indeed, if $a(j,n)=1$, then $a(j,n+1)>1$ by Eq.(\ref{frises_induction}); moreover, each $a(j,n)$ is a positive integer. We have, since $a(j,n_0)=a(j,n_0+p)$,
\begin{eqnarray*}
b(j)^2 &=& \left(\prod_{n_0\leq n <n_0+p} a(j,n)\right)\left( \prod_{n_0\leq n <n_0+p} a(j,n+1)\right) \\
&=& \prod_{n_0\leq n <n_0+p} a(j,n)a(j,n+1) \\
&=& \prod_{n_0\leq n <n_0+p} (1+(\prod_{j\rightarrow i} a(i,n)^{|C_{ij}|})(\prod_{i\rightarrow j} a(i,n+1)^{|C_{ij}|}))
\end{eqnarray*}
by Eq.(\ref{frises_induction}). Thus
\begin{eqnarray*}
b(j)^2 &>& \prod_{n_0\leq n <n_0+p}  ((\prod_{j\rightarrow i} a(i,n)^{|C_{ij}|})(\prod_{i\rightarrow j} a(i,n+1)^{|C_{ij}|})) \\
&=& (\prod_{j\rightarrow i} \prod_{n_0\leq n <n_0+p} a(i,n))^{|C_{ij}|}) (\prod_{i\rightarrow j}\prod_{n_0\leq n <n_0+p} a(i,n+1))^{|C_{ij}|}) \\
&=& (\prod_{j\rightarrow i} b(i)^{|C_{ij}|})  (\prod_{i\rightarrow j} b(i)^{|C_{ij}|}) \\
&=& \prod _{i\neq j} b(i)^{|C_{ij}|} .
\end{eqnarray*}
Taking logarithms, we obtain
$$
2 log(b(j)) > \sum _{i\neq j} log(b(i)) |C_{ij}|
$$
and we have a subadditive function which is not additive, since $b(j)>1$.

\n{\bf Case 2.} We assume now that some sequence $a(j,n)$ is unbounded. Then by 
Eq. (\ref{frises_induction}) and the connectedness of the underlying graph of the Cartan matrix, they are all unbounded.

We show, by using Lemma \ref{asymptotic}, that there exists an additive function for the Cartan matrix. We use freely the notations of this lemma. Define
$$b(j,n)=a(j,n)a(j,n+1)\cdots a(j,n+p-1).$$ 
Then 

\begin{equation}\label{bjpnk}
b(j,pn_k)\approx \lambda(j,0)^{n_k} n_k^{d(j,0)}\cdots \lambda(j,p-1)^{n_k} n_k^{d(j,p-1)}\approx \lambda(j)^{n_k} n_k^{d(j)}
\end{equation}

\n where $\lambda(j)=\lambda(j,0)\cdots\lambda(j,p-1)$ and $d(j)=d(j,0)+\cdots+d(j,p-1)$. Now
$$b(j,pn_k)^2=a(j,pn_k)a(j,pn_k+1)a(j,pn_k+1)a(j,pn_k+2)\cdots $$
$$ \cdots a(j,pn_k+p-1)a(j,pn_k).
$$
By the lemma, $a(j,pn_k)\approx a(j,pn_k+p)$. Thus
$$b(j,pn_k)^2 \approx\prod_{0\leq l<p}a(j,pn_k+l)a(j,pn_k+l+1).$$
Using Eq. (\ref{frises_induction}), we obtain
$$
b(j,pn_k)^2\approx \prod_{0\leq l<p} (1+(\prod_{j\rightarrow i} a(i,pn_k+l)^{|C_{ij}|})(\prod_{i\rightarrow j} a(i,pn_k+l+1)^{|C_{ij}|})).
$$
Let $$u_k=(\prod_{j\rightarrow i} a(i,pn_k+l)^{|C_{ij}|})(\prod_{i\rightarrow j} a(i,pn_k+l+1)^{|C_{ij}|}).$$ 
If $u_k$ is unbounded when $k\rightarrow \infty$, then by (i) in Lemma \ref{asymptotic}, there exists $i$ with: either $j\rightarrow i$, and $\lambda(i,l)>1$ or $d(i,l)\geq 1$; or $j\leftarrow i$, and $\lambda(i,l+1)>1$ or $d(i,l+1)\geq 1$. Then $\lim_{k\rightarrow\infty}u_k=\infty$ and $u_k\approx 1+u_k$. Otherwise, $u_k$ is bounded and by Lemma \ref{asymptotic}, $u_k$ is constant, therefore $u_k\approx 1+u_k$. Thus in both cases, $1+u_k\approx u_k$.

We deduce that 
$$
b(j,pn_k)^2 $$
$$\approx\prod_{0\leq l<p} (\prod_{j\rightarrow i} a(i,pn_k+l)^{|C_{ij}|})(\prod_{i\rightarrow j} a(i,pn_k+l+1)^{|C_{ij}|}) $$
$$\approx\prod_{0\leq l<p}(\prod_{j\rightarrow i}\lambda(i,l)^{n_k|C{ij}|} n_k^{d(i,l)|C_{ij}|}) (\prod_{i\rightarrow j} \lambda(i,l+1)^{n_k|C{ij}|} n_k^{d(i,l+1)|C_{ij}|}) $$
$$\approx(\prod_{j\rightarrow i}  \prod_{0\leq l<p} \lambda(i,l)^{n_k|C{ij}|} n_k^{d(i,l)|C_{ij}|}) (\prod_{i\rightarrow j}  \prod_{0\leq l<p} \lambda(i,l+1)^{n_k|C{ij}|} n_k^{d(i,l+1)|C_{ij}|}).$$
Since $d(i,0)=d(i,p)$ and $\lambda(i,0)=\lambda(i,p)$, we obtain
\begin{eqnarray*}
b(j,pn_k)^2 &\approx& (\prod_{j\rightarrow i} \lambda(i)^{n_k|C{ij}|} n_k^{d(i)|C_{ij}|}) (\prod_{i\rightarrow j} \lambda(i)^{n_k|C{ij}|} n_k^{d(i)|C_{ij}|}) \\
&\approx& \prod_{i\neq j}\lambda(i)^{n_k|C{ij}|} n_k^{d(i)|C_{ij}|}.
\end{eqnarray*}
Thus by Eq.(\ref{bjpnk}),
$$\lambda(j)^{2n_k}n_k^{2d(j)}\approx\prod_{i\neq j}\lambda(i)^{n_k|C{ij}|} n_k^{d(i)|C_{ij}|}.$$
Therefore, since $n_k$ tends to infinity with $k$, for $j=1,...,d$
$$\lambda(j)^2=\prod_{i\neq j}\lambda(i)^{|C{ij}|}
$$
and
$$
2d(j)=\sum_{i\neq j}d(i)|C_{ij}|.
$$
If the $d(j)$ are all positive, we have the additive function $d(j)$. If one of them is 0, then they are all 0, by connectedness of the graph and the above equation. In this case, $d(j,l)=0$ for any $j$ and $l$. Thus (ii) in the lemma ensures that for any $j$, some $\lambda(j,l)>1$ and therefore $\lambda(j)>1$. Taking logarithms, we find
$$
2log(\lambda(j))=\sum_{i\neq j}log(\lambda(i)))|C_{ij}|
$$
and we have therefore an additive function.

\begin{flushright}
$\square$
\end{flushright}

\section{Proof of Theorem 3}

\n Step 1. We prove first that the function $t$ given by Eq.(\ref{tilingformula}) is an $SL_2$-tiling of the plane. It is enough to show that for any $(u,v)\in \mathbf Z^2$, the determinant of the matrix $\left ( \begin{array}{ll} t(u,v) & t(u,v+1)\\ t(u+1,v) & t(u+1,v+1) \end{array} \right )$ is equal to 1.

By inspection of the figure below, where $k,l\ge 0$ and $w=x_1\cdots x_n$, $n\geq 0$ and $x_i\in\{x,y\}$,
$$
\begin{array}{cccccccccccc}
&&&&&&&&&&\bf 1&\bf 1 \\
&&&&&&&&&&.     &|      \\
&&&&&&&&&&.     &|     \\
&&&&&&&&&&.      &|    \\
&&&&&&&&&&\bf 1 &|  \\
&&&&&&&&&\bf 1&\bf 1&|  \\
&&&&&&&&.&&|&| \\
&&&&&&&w&&&|&|  \\
&&&&&&.&&&&|&|     \\
&&&&&\bf 1&&&&&|&|  \\
\bf 1&.&.&.&\bf 1&\bf 1&-&-&-&-&(u,v)&(u,v+1) \\
\bf 1&-&-&-&-&-&-&-&-&-&(u+1,v)&(u+1,v+1)
\end{array}
$$
\n it is seen that the words associated to the four points $(u,v)$, $(u,v+1)$, $(u+1,v)$ and $(u+1,v+1)$ are  respectively of the form $ywx$, $ywxy^{l}x$, $yx^{k}ywx$ and $yx^{k}ywxy^{l}x$.
Let $M=M(x_1\ldots x_n)$. Moreover, denote by $S(A)$ the sum of the coefficients of any matrix $A$. Then $t(u,v)=S(M)$, $t(u,v+1)=S(MM(x)M(y)^l)$, $t(u,v+1)=S(M(x)^kM(y)M)$ and moreover $t(u+1,v+1)=S(M(x)^kM(y)MM(x)M(y)^l)$. A straightforward computation, which uses the fact that $det(M)=1$, then shows that $t(u,v)t(u+1,v+1)-t(u,v+1)t(u+1,v)=1$.

\n Step 2. Clearly $t(u,v)> 0$ for any $(u,v)\in \mathbf Z^2$. Then it is easily deduced, by induction on the length of the word associated to $(u,v)$, that $t(u,v)$ is uniquely defined by the $SL_2$ condition. This proves that the tiling is unique. 

\begin{flushright}
$\square$
\end{flushright} 

\section{Properties of $SL_2$-tilings}

\subsection{ Rays and periodic frontiers}

Given a mapping $t: \mathbf Z^2\rightarrow R$, a point $M\in \mathbf Z^2$ and a nonzero vector $v\in \mathbf Z^2$, we consider the sequence $a_n=t(M+nv)$. Such a sequence will be called a {\em ray associated to} $t$. We call $M$ the {\em origin} of the ray and $v$ its {\em directing vector}. The ray is {\em horizontal} if $v=(1,0)$, {\em vertical} if $v=(0,-1)$ and {\em diagonal} if $v=(1,-1)$. 

We say that the frontier Eq.(\ref{frontier}) is {\em ultimately periodic} if for some $p\geq 1$, called a {\em period}, and some $n_0,n_0'\in \mathbf Z$, one has:

\n (i) for $n\geq n_0$, $x_n=x_{n+p}$;

\n (ii) for $n\leq n'_0$, $x_n=x_{n-p}$.

\begin{corollary} \label{rational}
If the frontier in Theorem 3 is ultimately periodic, then each ray associated to $t$, and whose directing vector is of the form $(a,b)$ with $ab\leq 0$,  is $\mathbf N$-rational.
\end{corollary}

\n\textit{Proof.}

Step 1. The points $M+nv$ are, for $n$ large enough, all above or all below the frontier, since the frontier is admissible and by the hypothesis on the directing vector. Since rationality is not affected by changing a finite number of values, we may, by symmetry, suppose that they are all below. 

Step 2. Let $w_n$ be the word associated to the point $M+nv$. 
By ultimate periodicity of the frontier, there exists an integer $q\geq 1$ and words $v_0,...,v_{q-1}, u'_0,...,u'_{q-1},u_0,...,u_{q-1}$ such that for any $i=0,...,q-1$ and for $n$ large enough, $w_{i+nq}={u'}_i^nv_iu_i^n$. 

Step 3. It follows from Theorem 3 that for some $2\times 2$ matrices $M'_i,N_i,M_i$ over $\mathbf N$, one has for any $i=0,\ldots,q-1$ and $n$ large enough, $a_{i+nq} = \lambda_i {M'_i}^nN_iM_i^n\gamma_i$, where $\lambda_i\in \mathbf N^{1\times 2},\gamma_i\in \mathbf N^{2\times1}$.

Step 4. Since rational series over $\mathbf N$ are closed under Hadamard product, each series $\sum_{n\in \mathbf N}a_{i+nq}x^n$ is rational over $\mathbf N$; indeed, such a series is by the formula in Step 3 an $\mathbf N$-linear combination of products of Hadamard products of two $\mathbf N$-rational series. Therefore $$\sum_{n\in \mathbf N}a_nx^n=\sum_{i=0,...,q-1} x^i (\sum_{n\in \mathbf N}a_{i+nq}(x^q)^n)$$
is  also $\mathbf N$-rational.
\begin{flushright}
$\square$
\end{flushright}

\subsection{Symmetric frontiers, perfect squares and quadratic relations}

Given a finite or infinite word $w$ on the alphabet $\{x,y\}$, we call {\em transpose} of $w$, and denote it by $^tw$ the word obtained by reversing it and exchanging $x$ and $y$. For instance, $^t(xyyxy)=xyxxy$. If $w$ is a right infinite word, then its transpose is a left infinite word.

\begin{lemma}
Consider an admissible frontier of the form ${}^tsyx^hys$, embedded in the plane, with $h\in\mathbf N$. Let $I,J,K$ be the points of the plane defined as follows: $I$ corresponds to the point between $x^h$ and $y$ on the frontier; $J$ (resp. $K$) is immediately below $I$ (resp. $J$); see the figure. Let $i_n$ (resp. $j_n$, $k_n$) be the horizontal (resp. diagonal) ray of origin $I$ (resp. $J$, $K$) of the $SL_2$-tiling corresponding to the frontier. Then for any $n\in \mathbf N$,
$$j_n=(h+1)i_n^2   \,\,\mbox{and}  \,\, k_{n}+1=(h+1)i_ni_{n+1}.$$
\end{lemma}

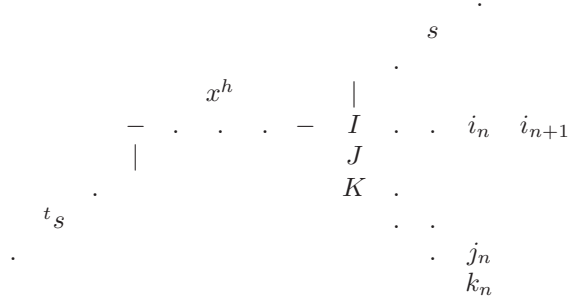
$$
\begin{array}{cccccccccccccccc}
&&&&&&&&&&&. \\
&&&&&&&&&&s \\
&&&&&&&&&. \\
&&&&&x^h&&&| \\
&&&-&.&.&.&-&I &.&.&i_n&i_{n+1}\\
&&&|&&&&&J \\
&&.&&&&&&K&. \\
&{}^ts&&&&&&&&.&.\\
. &&&&&&&&&&.&j_n\\
&&&&&&&&&&&k_n
\end{array}
$$

\b
\n\textit{Proof.}
Step 1. We have 
$$M(yx^hy)=\left(\begin{array}{cc}1&0\\1&1\end{array}\right) \left(\begin{array}{cc}1&h\\0&1\end{array}\right) \left(\begin{array}{cc}1&0\\1&1\end{array}\right)=\left(\begin{array}{cc}h+1&h\\h+2&h+1\end{array}\right).$$
The quadratic form associated to this matrix (which is not symmetric) is therefore
$$(a,b) \left(\begin{array}{cc}h+1&h\\h+2&h+1\end{array}\right) \left(\begin{array}{c}a\\b\end{array}\right)=(h+1)a^2+(2h+2)ab+(h+1)b^2=(h+1)(a+b)^2.$$

Step 2. The words associated to $i_n$, $j_n$, $i_{n+1}$, $k_n$ are respectively of the form $yvx$, $y\,{}^tvyx^hyvx$, $yvxy^kx$, $yx^ky\,{}^tvyx^hyvx$; see the figure.

$$
\begin{array}{cccccccccccccccccccc}
&&&&&&&&&&&&&&&&&\bf 1&\bf 1 \\
&&&&&&&&&&&&&&&&&. \\
&&&&&&&&&&&&&&&&k&. \\
&&&&&&&&&&&&&&&&&. \\
&&&&&&&&&&&&&&&&\bf 1&\bf 1 \\
&&&&&&&&&&&&&&&. \\
&&&&&&&&&&&&&&v \\
&&&&&&&&&&&&&. \\
&&&&&&&&&&h&&\bf 1 \\
&&&&&&&&\bf 1 &.&.&.&\bf 1&&&&&i_n&i_{n+1}\\
&&&&&&&&\bf 1 \\
&&&&&&&. \\
&&&&&&{}^tv \\
&&&&&. \\
&&k&& \bf 1\\
\bf 1&.&.&.&\bf 1 &&&&&&&&&&&&&j_n\\
\bf 1 &&&&&&&&&&&&&&&&&k_n
\end{array}
$$

Step 3. Thus $i_n=(1,1)M(v) \left(\begin{array}{c}1\\1\end{array}\right)$ and $j_n=(1,1)M({}^tvyx^hyv) \left(\begin{array}{c}1\\1\end{array}\right)$. Let $\left(\begin{array}{c}a\\b\end{array}\right)=M(v)\left(\begin{array}{c}1\\1\end{array}\right)$. Then $i_n=a+b$ and $j_n=(1,1)M({}^tv)M(yx^hy)M(v)\left(\begin{array}{c}1\\1\end{array}\right)=(a,b)M(yx^hy)\left(\begin{array}{c}a\\b\end{array}\right)$, since $M({}^tv)={}^tM(v)$. Thus by Step 1, $j_n=(h+1)(a+b)^2=(h+1)i_n^2$.

Step 4. Let $M(v)=\left(\begin{array}{cc}p&q\\r&s\end{array}\right)$. Then 
$$
\begin{array}{lll}
k_n&=&(1,1)M(x^ky\,{}^tvyx^hyv)\left(\begin{array}{c}1\\1\end{array}\right)\\
&=&\left(\begin{array}{cc}1&k\\0&1\end{array}\right) \left(\begin{array}{cc}1&0\\1&1\end{array}\right) \left(\begin{array}{cc}p&r\\q&s\end{array}\right) \left(\begin{array}{cc}h+1&h\\h+2&h+1\end{array}\right) \left(\begin{array}{cc}p&q\\r&s\end{array}\right).
\end{array}
$$
Furthermore, $i_n=p+q+r+s$ and 
$$
\begin{array}{lll}
i_{n+1}&=&(1,1)M(vxy^k)\left(\begin{array}{c}1\\1\end{array}\right)\\
&=&(1,1)\left(\begin{array}{cc}p&q\\r&s\end{array}\right) \left(\begin{array}{cc}1&1\\0&1\end{array}\right) \left(\begin{array}{cc}1&0\\k&1\end{array}\right) \left(\begin{array}{c}1\\1\end{array}\right).
\end{array}
$$
A straightforward computation then shows that $(h+1)i_ni_{n+1}-k_n=ps-rq$. Since $det(M(v))=1$, the lemma is proved.

\begin{flushright}
$\square$
\end{flushright}

\n {\bf Remark} Denote by $k'_n$ the value of the tiling in the point immediately to the right of the value $j_n$. It is easily shown that one has also $k'_{n}-1=(h+1)i_ni_{n+1}$. Therefore, the $2\times 2$ matrix $\left(\begin{array}{cc}j_n&k'_n\\k_n&j_{n+1}\end{array}\right)$, which appears as a connected submatrix of the tiling, encodes a pythagorean triple: indeed, $(j_{n+1}+j_n, j_{n+1}-j_n, k_n+k'_n)$ is such a triple, because $(j_{n+1}-j_n)^2+(k_n+k'_n)^2=(h+1)^2(i_{n+1}^2-i_n^2)^2+(h+1)^2(2i_ni_{n+1})^2=(h+1)^2(i_{n+1}^2+i_n^2)^2=(j_{n+1}+j_n)^2$. See for example the tiling given in Section 3 and its submatrices $\left(\begin{array}{cc}1&3\\1&2^2\end{array}\right)$, representing the triple $(5,3,4)$, or $\left(\begin{array}{cc}2^2&11\\9&5^2\end{array}\right)$, representing the triple $(29,21,20)$ (here $h=0$).

\begin{lemma}\label{quadratic}
Let $h,h'\in\mathbf N$. Consider a frontier of the form $f=s'xy^{h'}xwyx^hys$, where $w\in\{x,y\}^*$, such that ${}^ts=s'xy^{h'}xw$ and ${}^ts'=wxyx^hys$. Let $P_0,...,P_k$ be the points corresponding to $w$, with $k$ the length of $w$, and $b(j,n)$ be the diagonal rays of origin $P_j$, for $j=0,...,k$. Let $I,J,K,I',J',K'$ be the points defined as follows: $I$ (resp. $I'$) corresponds to the point on the frontier between $x^h$ and $y$ (resp. $x$ and $y^{h'}$), $J$ (resp. $K$) is immediately below $I$ (resp. $J$), $J'$ (resp. $K'$) is immediately to the right of $I'$ (resp. of $J'$). Let $i_n$ (resp. $i'_n$, resp. $ j_n,k_n,j'_n,k'_n$) denote the horizontal (resp. vertical, resp. diagonal) ray of origin $I$ (resp. $I'$, resp. $J,K,J',K'$). Then  $j_n=(h+1)i_n^2$, $k_n+1=(h+1)i_ni_{n+1}$, $j'_n=(h'+1){i'}_n^2, k'_n+1=(h'+1)i'_ni'_{n+1}$. Moreover, for $j=1,\ldots ,k-1$, one has, with $w=x_1\cdots x_k$: for any $n\in \mathbf N$, $b(j,n)b(j,n+1)=1+B$, where 
$$
B= \left\{ 
\begin{array}{lll} b(j-1,n+1)b(j+1,n)&{if} &x_jx_{j+1}=xx \\
b(j-1,n+1)b(j+1,n+1)&{if} & x_jx_{j+1}=xy  \\
b(j-1,n)b(j+1,n)&{if} & x_jx_{j+1}=yx\\
b(j-1,n)b(j+1,n+1)&{if} & x_jx_{j+1}=yy
\end{array}
\right.
$$
\end{lemma}

Note that the points given in the lemma need not be all distinct. The lemma is illustrated by the following figure, where $h=1$, $h'=0$, and $w=xyxx$; in bold are represented the 1's corresponding to the factors $xy^{h'}x=xx$ and $yx^hy=yxy$ of the frontier.

$$
\begin{array}{lllllllllllllllllllllllllllllllllllllllllllllll}
&&&&&&&&&1&1 \\
&&&&&&&&&1&2 \\
&&&\cdots&&&&&&1&3 \\
&&&&&&&&1&1&4 \\
&&&&&&&&1&2&9 \\
&&&&&&&&\bf 1&3&14 \\
&&&&&&&\bf 1&{\bf 1}(I)&4&19 \\
&&&&&1(P_2)&1(P_3)&{\bf 1}(P_4)&2\cdot1^2(J)&9&43 \\
&&\bf 1&{\bf 1}(I')&{\bf 1^2}(J'=P_0)&1(K'=P_1)&2&3&7(K)&2\cdot4^2&153 \\
&1&1&2&3&2^2&9&14&33&151&2\cdot19^2\\
&1&2&5&8&11&5^2&39\\
1&1&3&8\\
&&&&&&&&\cdots
\end{array}
$$

\b
\n\textit{Proof.} This follows from Lemma 1, and its symmetrical statement, together with the inspection of the figure below, which shows the four different possible configurations.
$$
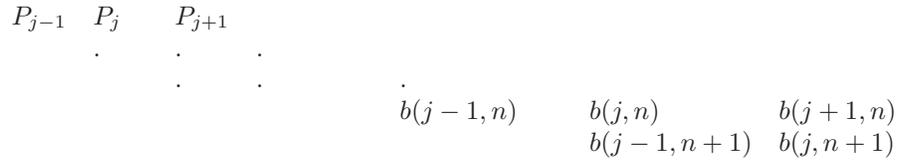
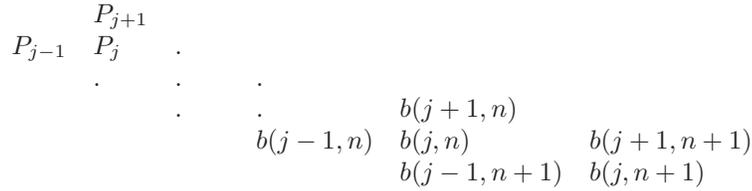
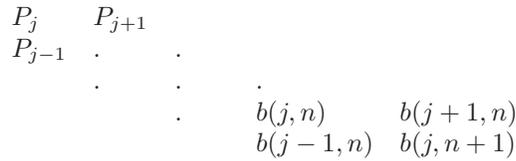
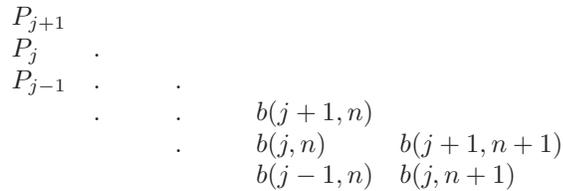
\begin{array}{lllllllllllllll}
P_{j-1}&P_j&P_{j+1} \\
&.&.&. \\
&&.&.&. \\
&&&&b(j-1,n)&b(j,n)&b(j+1,n) \\
&&&&&b(j-1,n+1)&b(j,n+1)\\
\\
\\
&P_{j+1} \\
P_{j-1}&P_j&. \\
&.&.&. \\
&&.&.&b(j+1,n) \\
&&&b(j-1,n)&b(j,n)&b(j+1,n+1) \\
&&&&b(j-1,n+1)&b(j,n+1)
\\
\\
P_j&P_{j+1} \\
P_{j-1}&.&. \\
&.&.&. \\
&&.&b(j,n)&b(j+1,n) \\
&&&b(j-1,n)&b(j,n+1) \\
\\
\\
P_{j+1} \\
P_j &.\\
P_{j-1}&.&. \\
&.&.&b(j+1,n) \\
&&.&b(j,n)&b(j+1,n+1) \\
&&&b(j-1,n)&b(j,n+1) \\
\end{array}
$$

\begin{flushright}
$\square$
\end{flushright}

\begin{lemma}\label{periodic-exists}
Let $h,h'$ in $\mathbf N$ and $w$ in $\{x,y\}^*$. Then there exists a periodic frontier $f=s'xy^{h'}xwyx^hys$ satisfying the hypothesis of Lemma \ref{quadratic}.
\end{lemma}

\b
\n\textit{Proof.}
Let indeed
\begin{equation}\label{periodic_frontier}
f={}^\infty(wxy^hx\,{}^twxy^{h'}x)(wyx^hy\,{}^twyx^{h'}y)^{\infty}.
\end{equation}
In other words, we take $s=({}^twyx^{h'}ywyx^hy)^\infty$ and $s'={}^\infty(xy^{h'}xwxy^hx\,{}^tw)$. Then ${}^ts={}^\infty(xy^hx\,{}^twxy^{h'}xw)={}^\infty(xy^{h'}xwxy^hx\,{}^tw)xy^{h'}xw=s'xy^{h'}xw$. Similarly ${}^ts'=wyx^hys$.
\begin{flushright}
$\square$
\end{flushright}

\section{Proof of Theorem 2}

We completely omit the case of Dynkin diagrams, since Theorem 2 in this case follows immediately from the finiteness of the set of cluster variables, see \cite{FZ2}. 

\subsection{The case $\tilde A_m$}

Let $1,...,m+1$ be the vertices of the graph $\tilde A_m$, with edges $\{j,j+1\}$, $j=1,...,m+1$, with $j+1$ taken $mod. m+1$. An acyclic orientation being given, let $x_j=x$ if the orientation is $j\rightarrow j+1$ and $x_j=y$ if it is $j\leftarrow j+1$. Let $a(j,n)$ be the sequences of the frise, $j=1,\ldots,m+1$. We extend the notation $a(j,n)$ to $j\in \mathbf Z$ by taking $j$ $mod. m+1$. Let $w$ be the word $x_1\cdots x_{m+1}$, which encodes the orientation. Then $^\infty w^\infty$ is an admissible frontier; indeed, $x$ and $y$ appear both in $w$, since the orientation is acyclic. Embed this frontier into the plane and denote by $P_j$ , with $j\in \mathbf Z$, the successive points of this embedding, in such a way that  $P_j$ corresponds to the point between $x_{j-1}$ and $x_j$, with $j$ taken $mod.\, m+1$. Let $t$ be the tiling given by Theorem 3. Because of the periodicity of the frontier, the diagonal ray $b(j,n)$ of origin $P_j$, $j\in \mathbf Z$, depends only on the class of $j$ $mod. \, m$.

We claim that the ray $b(j,n)$ is equal to $a(j,n)$, for $j=1,..,m+1$. This is true for $n=0$, since both are equal to 1. It is enough to show that $b(j,n)$ satisfies Eq.(1). Fix $j=1,...,m+1$. We have four cases according to the relative positions of $P_{j-1},P_j,P_{j+1}$ (see the figure in the proof of Lemma \ref{quadratic}).

They correspond to the four possible values of the couple $(x_{j-1},x_j)$: $$ (x,x), (x,y), (y,x),(y,y).$$
By definition of $w$, these four cases correspond to the four possible orientations:
$$j-1\rightarrow j\rightarrow j+1, \,j-1\rightarrow j\leftarrow j+1,$$
$$\, j-1\leftarrow j\rightarrow j+1,\, j-1\leftarrow j\leftarrow j+1.
$$
Thus, by Eq.(1), they correspond to the four induction formulas $a(j,n+1)=\frac{1+A}{a(j,n)}$, where $A$ takes one of the four possible values:
$$
a(j-1,n+1)a(j+1,n), \quad a(j-1,n+1)a(j+1,n+1), $$
$$a(j-1,n)a(j+1,n), \quad a(j-1,n)a(j+1,n+1).
$$
Regarding the tiling, these four cases correspond to the four possible configurations, shown in the same figure.
Hence, by the $SL_2$-condition, they correspond to the four induction formulas for $b(j,n)$: $b(j,n+1)=\frac{1+B}{b(j,n)}$, where $B$ takes one of the four possible values:
$$
 b(j-1,n+1)b(j+1,n), b(j-1,n+1)b(j+1,n+1),$$
$$ b(j-1,n)b(j+1,n), b(j-1,n)b(j+1,n+1).
$$
This concludes the proof, by using Corollary 1. \begin{flushright}
$\square$
\end{flushright}

The proof is illustrated in the tiling below and in Figure \ref{friseA3tilde}, for a specific orientation of $\tilde A_3$.

$$\begin{array}{lllllllllllllllllllllllllllllll}
&&&&&&&&&&&&1 \\
&&&&&\ldots&&&&1&1&1&1 \\
&&&&&&1&1&1&1&2&3&4 \\
&&&1&1&1&1&2&3&4&9&14&19 \\
1&1&1&1&2&3&4&9&14&19&43&67 &&\ldots \\
1&2&3&4&9&14&19&43&67 \\
&&&&\ldots
\end{array}
$$

\begin{figure}[htbp] 
\begin{center}

 \scalebox{0.6}{\input{friseA3tilde.pstex_t}}

\caption{Quiver $\tilde A_3$ and frise}
\label{friseA3tilde}
\end{center}
\end{figure}

\subsection{The case $\tilde D_m$}

\begin{figure}[htbp]
\begin{center}

 \scalebox{.6}{\input{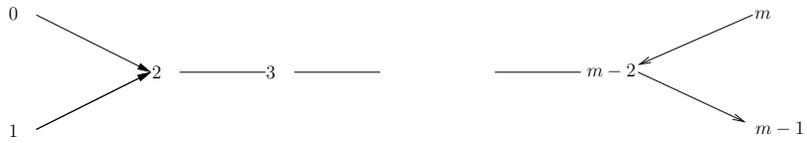}}

\caption{Quiver $\tilde D_m$}
\label{Dmtilde}
\end{center}
\end{figure}

We consider an orientation of $\tilde D_m$, of the form shown in Figure \ref{Dmtilde}, where the orientations of the edges $\{i,i+1\}$ for $i=2,...,m-3$ are arbitrary. The other cases, which differ from the case considered here by changing the orientation of the forks, are similar.

\begin{figure}[htbp]
\begin{center}

 \scalebox{.7}{\input{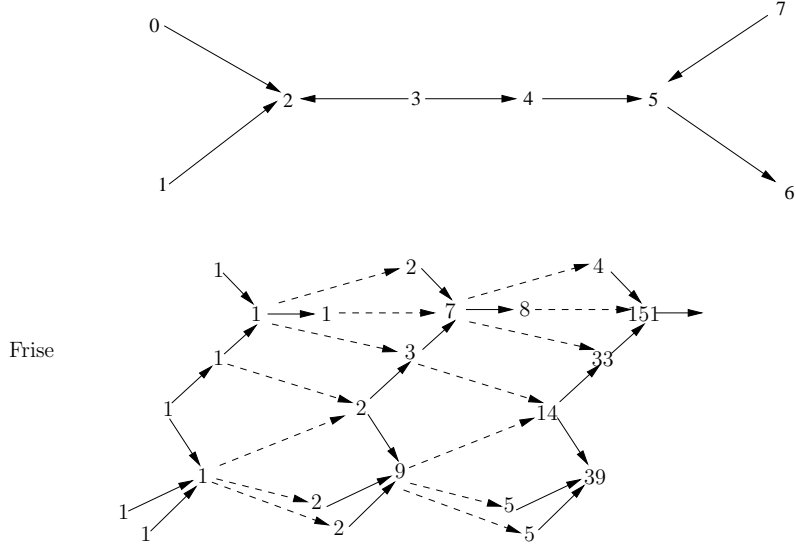}}

\caption{Quiver $\tilde D_7$}
\label{friseD7tilde}
\end{center}
\end{figure}

An example is shown in  Figure \ref{friseD7tilde}. The reader may recognize that this frise is encoded in the tiling shown in Section 6.2 after the statement of Lemma 3.

\n Step 1. We define a word $w$, which  encodes the orientation, as follows: $w=x_1x_2\cdots x_{m-3}$, with
$$
x_i=\left\{ 
\begin{array}{ll}
x & \mbox{if $i\rightarrow i+1$} \\
y & \mbox{if $i\leftarrow i+1$}.
\end{array}
\right.
$$

Note that by our choice in Figure \ref{Dmtilde}, $x_1=x$. In the example, $w=xyxx$. We consider now the frontier given by Lemma \ref{periodic-exists}, with $h'=0$ and $h=1$; see Eq.(\ref{periodic_frontier}). The associated rays $b(j,n)$ are all rational by Corollary \ref{rational}.

\n Step 2. Taking the notations of Lemma \ref{quadratic}, we have for any $n\in\mathbf N$:
$$
j'_n=b(0,n), \, k'_n=b(1,n), \, k_n=b(m-3,n+1). 
$$
For $j'_n$ and $k'_n$ this follows from $h'=0$ and therefore $J'=P_0$, $K'=P_1$, hence $j'_n$ (resp. $k'_n$) and $b(0,n)$ (resp. $b(1,n)$) are diagonal rays with the same origin. Moreover $w$ is of length $k$ in Lemma \ref{quadratic} and here of length $m-3$, hence $k=m-3$. Since $h=1$, we have the configuration shown below.

$$
\begin{array}{lllll}
&1 \\
1&1(I) \\
1(P_{m-3})&(J) \\
&(K)

\end{array}
$$

This implies $k_n=b(m-3,n+1)$.

\n Step 3. In accordance with Lemma \ref{quadratic}, we have for $j=1,\ldots, k$,
$$
x_j=\left\{ 
\begin{array}{ll}
x & \mbox{if $[P_{j-1},P_j]$ is horizontal} \\
y & \mbox{if $[P_{j-1},P_j]$ is vertical}.
\end{array}
\right.
$$

\n Step 4.
We define $m+1$ sequences $a'(j,n)$ for $j=0,\ldots m$. First, for any $n\in\mathbf N$,
$$
a'(0,n)=a'(0,n)=a'(1,n)=i'_n.
$$
Now, for $j=2,...,m-2$,
$$
a'(j,n)=b(j-1,n).
$$
Furthermore, 
$$
a'(m-1,n)=\left\{ 
\begin{array}{ll}
i_n & \mbox{if $n$ is even} \\
2i_n & \mbox{if $n$ is odd}.
\end{array}
\right.
$$
Finally, for $n\ge 1$,
$$
a'(m,n)=\left\{ 
\begin{array}{ll}
i_{n-1} & \mbox{if $n$ is even} \\
2i_{n-1} & \mbox{if $n$ is odd},
\end{array}
\right.
$$
with $a'(m,0)=1$.

Observe that for any $n\in\mathbf N$, $a'(m-1,n)a'(m,n+1)=2i_n^2$. Moreover, since the sequences $i_n$, $i'_n$ and $b(j,n)$ are rays, they are rational. Hence, so are the sequences $a'(j,n)$ (for $a'(m-1,n)$ and $a'(m,n)$, this follows from standard constructions on rational sequences).

\n Step 5. In order to end the proof, it is enough to show that $a(j,n)=a'(j,n)$. First note that $a(j,0)=a'(j,0)$, as is easily verified. Thus it suffices to show that $a'(j,n)$ satisfies the same recursion formula as $a(j,n)$, that is, Eq.(\ref{frises_induction}).

\n Step 6. We have
$$
\begin{array}{llll}
a(0,n+1)a(0,n)&=& i'_{n+1}i'_n & \mbox{by Step 4} \\
&=& 1+k'_n & \mbox{by Lemma \ref{quadratic}} \\
&=& 1+b(1,n) & \mbox{by Step 2} \\
&=& 1+a'(2,n) & \mbox{by Step 4}.
\end{array}
$$
Similarly
$$a'(1,n+1)a'(1,n)=1+a'(2,n).$$
This is the good recursion for $a'(0,n)$ and $a'(1,n)$ since, by Figure \ref{Dmtilde} and Eq.(1), 
$$a(0,n+1)a(0,n)=1+a(2,n),$$
$$a(1,n+1)a(1,n)=1+a(2,n).$$

\n Step 7. Let
$$
n'=\left\{ 
\begin{array}{ll}
n & \mbox{if $x_2=x$} \\
n+1& \mbox{if $x_2=y$}.
\end{array}
\right.
$$
We have by Step 4, $a'(2,n+1)a'(2,n)=b(1,n+1)b(1,n)$. Looking at the figures below, where we use Step 3:

\n Case $n'=n,x_2=x$
$$
\begin{array}{llllll}
P_0&P_1&P_2&b(1,n)&b(2,n)\\
&&&b(0,n+1)&b(1,n+1)
\end{array}
$$
Case $n'=n+1, x_2=y$
$$
\begin{array}{lllll}
&P_2&b(1,n)&b(2,n+1)\\
P_0&P_1&b(0,n+1)&b(1,n+1)
\end{array}
$$
we see that this is equal to $1+b(0,n+1)b(2,n')$. Using Step 2, Step 4, Lemma \ref{quadratic} then Step 4 again, we obtain
$$
\begin{array}{lll}
b(0,n+1)b(2,n')&=&j'_{n+1}a'(3,n')\\
&=&{i'}_{n+1}^2a'(3,n')\\
&=&a'(0,n+1)a'(1,n+1)a'(3,n').
\end{array}
$$
Thus $a'(2,n+1)a'(2,n)=1+a'(0,n+1)a'(1,n+1)a'(3,n')$. This corresponds to Eq.(1) for $(2,n)$, that is, $a(2,n+1)a(2,n)=1+a(0,n+1)a(1,n+1)a(3,n')$.

\begin{figure}[htbp]
\begin{center}

 \scalebox{0.7}{\input{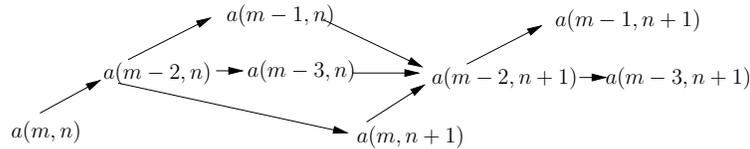}}

\caption{Frise in the case $n''=n$}
\label{FriseA}
\end{center}
\end{figure}

\begin{figure}[htbp]
\begin{center}
 
 \scalebox{0.7}{\input{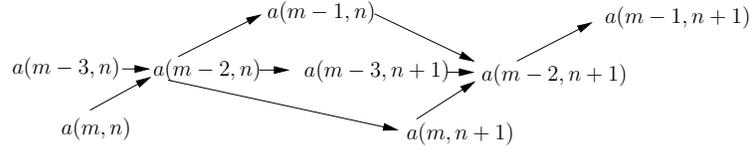}}

\caption{Frise in the case $n''=n+1$}
\label{FriseB}
\end{center}
\end{figure}

\n Step 8. Let
$$
n''=\left\{ 
\begin{array}{ll}
n & \mbox{if $x_{m-3}=y$.} \\
n+1& \mbox{if $x_{m-3}=x$}.
\end{array}
\right.
$$
We have 
$$
\begin{array}{llll}
a'(m-2,n+1)a'(m-2,n)&=&b(m-3,n+1)b(m-3,n)&\mbox{by Step 4}\\
&=&1+b(m-4,n'')j_n,
\end{array}
$$
where the second equality follows from figures \ref{FriseA} and \ref{FriseB}, which use Step 3:

\n Case $n''=n, x_{m-3}=y$
$$
\begin{array}{lllll}
P_{m-3}&J&&b(m-3,n)&j_n\\
P_{m-4}&&&b(m-4,n)&b(m-3,n+1)
\end{array}
$$
\n Case $n''=n+1, x_{m-3}=x$
$$
\begin{array}{llllll}
P_{m-4}&P_{m-3}&J&&b(m-3,n)&j_n\\
&&&&b(m-4,n+1)&b(m-3,n+1)
\end{array}
$$

By Lemma \ref{quadratic}, $j_n=2i_n^2$ and by Step 4, $2i_n^2=a'(m-1,n)a'(m,n+1)$. Furthermore, $b(m-4,n'')=a'(m-3,n'')$. Thus
$$
a'(m-2,n+1)a'(m-2,n)=1+a'(m-3,n'')a'(m-1,n)a'(m,n+1).
$$
This equality corresponds to Eq.(1) for $a(m-2,n)$, if we look at the two figures \ref{FriseA} and \ref{FriseB}, representing the frise in the two cases $n''=n, x_{n-3}=y$ and $n''=n+1, x_{n-3}=x$.

\n Step 9. We have 
$$
\begin{array}{lllll}
a'(m-1,n+1)a,(m-1,n)&=&2i_ni_{n+1}&\mbox{by Step 4}\\
&=&1+k_n&\mbox{by Lemma \ref{quadratic}}\\
&=&1+b(m-3,n+1)&\mbox{by Step 2}\\
&=&1+a'(m-2,n+1)&\mbox{by Step 4}\\
\end{array}
$$
in accordance with Eq.(1), which gives $a(m-1,n+1)a(m-1,n)=1+a(m-2,n+1)$.
Moreover, for $n\geq 1$, 
$$
\begin{array}{lllll}
a'(m,n+1)a'(m,n)&=&2i_ni_{n-1}&\mbox{by Step 4}\\
&=&1+k_{n-1}&\mbox{by Lemma \ref{quadratic}}\\
&=&1+b(m-3,n)&\mbox{by Step 2}\\
&=&1+a'(m-2,n)&\mbox{by Step 4},\\
\end{array}
$$
in accordance with Eq.(1), which gives $a(m,n+1)a(m,n)=1+a(m-2,n)$, noting that for $n=0$:
$a'(m,1)a'(m,0)=2i_0$ (by Step 4)$=2$ and $1+a'(m-2,0)=1+1=2$.

\subsection{Cases $\tilde B_m$, $\tilde C_m$, $\tilde BC_m$, $\tilde BD_m$, $\tilde CD_m$.}

Each frise in these cases is reduced to a frise of type $\tilde A$ or $\tilde D$. Precisely: $\tilde B_m$ is reduced to $\tilde D_{m+2}$; $\tilde C_m$ is reduced to $\tilde A_{2m}$; $\tilde {BC}_m$ is reduced to $\tilde D_{2m+2}$; $\tilde {BD}_m$ is reduced to $\tilde D_{m+1}$; $\tilde {CD}_m$ is reduced to $\tilde D_{2m}$.

\begin{figure}[htbp]
\begin{center}

 \scalebox{.6}{\input{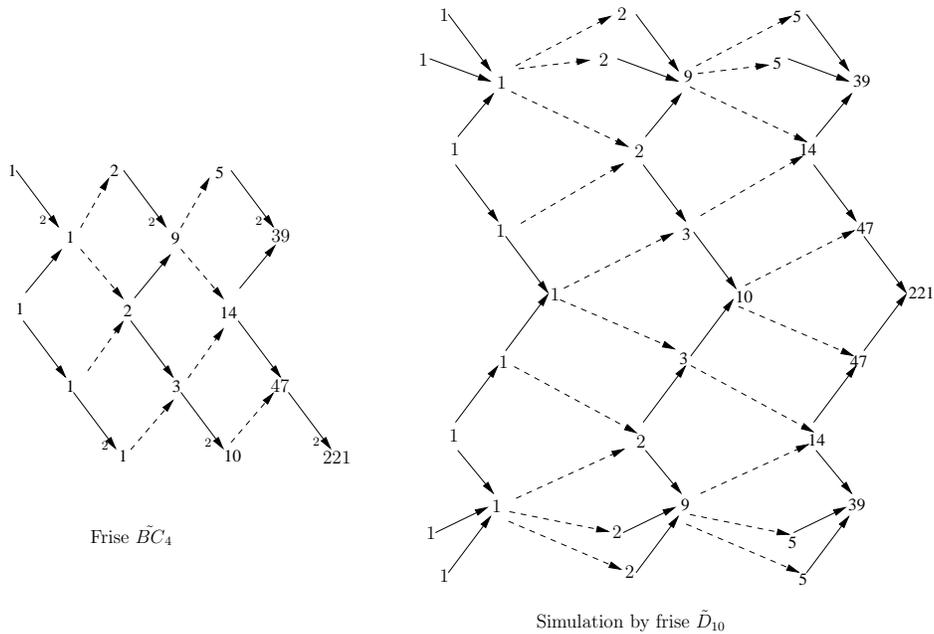}}

\caption{Frise simulation}
\label{BC4tilde}
\end{center}
\end{figure}

We give no formal proof, but an example; it should convince the reader. In this example, we show how a frise of type $\tilde {BC}_4$ can be simulated by a frise of type $\tilde D_{10}$. See Figure \ref{BC4tilde}.

\section{Frises and tilings with variables}

As mentioned in the introduction, if we replace the inital values $a(j,0)$ in
Section 3 by commuting variables, and keep the recurrence of Eq.(1) unchanged,
we obtain frises of variables.  The variables therefore obtained are
usual cluster variables, in the sense of Fomin and Zelevinsky.  It is
well-known that all cluster variables, but finitely many of them, can be
obtained in this way (those not obtained in this way being the cluster
variables corresponding to the exceptional objects, in the corresponding
cluster category, lying in tubes or $\mathbb{Z}A_{\infty}$ components).

Likewise, we generalize the $SL_2$-tilings; these are simply fillings of the
discrete plane by elements of a ring $R$ such that each $2\times 2$ connected
minor is of determinant $1$. We generalize Theorem 3 by putting variables
on the frontier.

\subsection{Case $\tilde A_m$}

We call ({\em generalized}) {\em frontier} a bi-infinite sequence 
\begin{equation} \label{frontiervariable}
\ldots x_{-2}a_{-2}x_{-1}a_{-1}x_0a_0x_1a_1x_2a_2x_3a_3 \ldots
\end{equation}
where $x_i\in \{x,y\}$ and $a_i$ are variables, for any $i\in \mathbf Z$. It is called {\em admissible} if there are arbitrarily large and arbitrarily small $i$'s such that $x_i=x$, and similarly for $y$; in other words, none of the two sequences $(x_n)_{n\geq0}$ and $(x_n)_{n\leq0}$ is ultimately constant. The $a_i$'s are called the {\em variables} of the frontier. Each frontier may be embedded into the plane: the variables label points in the plane, and the $x$ (resp. $y$) determine a bi-infinite discrete path, in such a way that $x$ (resp. $y$) corresponds to a segment of the form $[(a,b),(a+1,b)]$ (resp $[(a,b),(a,b+1)]$). 
For example, corresponding  to the frontier $\ldots a_{-2}xa_{-1}xa_0ya_1ya_2xa_3 \ldots$ is given below:

$$
\begin{array}{lllllllllllllll} 
&&&&&&&&&. \\
&&&&&&&&.    \\
&&&&&&&.       \\
&&&&&a_2&a_3   \\
&&&&&a_1  \\
&&&a_{-2}&a_{-1}&a_0  \\
&&. \\
&. \\
. \\

\end{array}
$$
Formally we do as follows. We define a partial function$f$ from $\mathbf Z^2$ into the semiring of Laurent polynomials over $\mathbf N$ generated by the variable, defined up to translation, as follows: fix some $(k,l)\in \mathbf Z^2$ and $i \in \mathbf Z$; then $f(k,l)=a_i$; moreover, if $x_{i+1}x_{i+2}...x_p$ labels the discrete path from $(k,l)$ to $(k',l')$, then $f(k',l')=a_p$; furthermore, if $x_px_{p+1}...x_i$ labels the path from $(k',l')$ to $(k,l)$, then $f(k',l')=a_{p-1}$.
  
We see below that an admissible frontier, embedded into the plane, may be extended to an $SL_2$-tiling. For this, we need the following notation. Let 
$$M(a,x,b)=\left(\begin{array}{cc}a&1\\0&b\end{array}\right) \,\mbox{and} \, \,M(a,y,b)=\left(\begin{array}{cc}b&0\\1&a\end{array}\right).$$

Note that these matrices reduce to the matrices $M(x)$ and $M(y)$ when the variables $a$ and $b$ are set to 1.
 
Given an admissible frontier, embedded in the plane as explained previously, let $(u,v)\in \mathbf Z^2$. Then we obtain a finite word, which is a factor of the frontier, by projecting the point $(u,v)$ horizontally and vertically onto the frontier. We call this word the {\em word} of $(u,v)$. It is illustrated in the figure below, where the word of the point $M=(k,l)$ is $a_{-3}ya_{-2}ya_{-1}ya_0xa_1xa_2ya_3xa_4$:
$$
 \begin{array}{ccccccccccccccccc} 
&&&&&&&&&&&&. \\
&&&&&&&&&&&.    \\
&&&&&&&&&&.       \\
&&&&&&&a_3&a_4&a_5 \\
&&&&&a_0&a_1&a_2  &|\\
&&&&&a_{-1}  &&&|\\
&&&&&a_{-2}  &&&|\\
&&&&a_{-4}&a_{-3} &-&-&M\\
&&. \\
&. \\
\end{array}
$$
We define the word of a point only for points below the frontier; for points above, things are symmetric and we omit this case. We call {\em denominator} of the point $M$ the product of the variables of its word, excluding the two extreme ones. In the example, its denominator is $a_{-2}a_{-1}a_0a_1a_2a_3$.

\begin{theorem}\label{tiling-variables}
Given an admissible frontier, there exists a unique $SL_2$-tiling  $t$ of the plane, with values in the semiring of Laurent polynomials over $\mathbf N$ generated by the variables lying on the frontier, extending the embedding of the frontier into the plane. It is defined, for any point $(u,v)$ below the frontier, with associated word $a_0x_1a_1x_2...x_{n+1}a_{n+1}$, where $n\geq 1$ and $x_i\in\{x,y\}$, by the formula
\begin{equation}\label{tilingformulvariable}
t(u,v)=\frac{1}{a_1a_2...a_n} (1,a_0)M(a_1,x_2,a_2)M(a_2,x_3,a_3)\cdots M(a_{n-1},x_n,a_n)\left(\begin{array}{c}1\\a_{n+1}\end{array}\right).
\end{equation}
\end{theorem}

In order to prove the theorem, we need two lemmas, where $R$ denotes some commutative ring. We extend the notation $M(a,x,b)$ and $M(a,y,b)$ for $a,b$ in $R$.

\begin{lemma} 
\label{lemmapqrsa}
\n (i) Let $A\in R^{2\times 2}$, $\lambda, \lambda'\in R^{1\times 2}$, $\gamma, \gamma'\in R^{2\times 1}$, and define $p=\lambda A\gamma$, $q=\lambda A\gamma'$, $r=\lambda' A\gamma$,  $s=\lambda' A\gamma'$. Then $det\left(\begin{array}{cc}p&q\\r&s\end{array}\right) =det(A)det\left(\begin{array}{c}\lambda\\\lambda'\end{array}\right) det(\gamma,\gamma')$.

\n (ii) Let $a,b_1,...,b_k,b\in R$, $\lambda'=(1,a)M(b_1,x,b_2)\cdots M(b_{k-1},x,b_k)M(b_k,y,b)$ and $\lambda=(1,b_k)$. Then $det\left(\begin{array}{c}\lambda'\\\lambda\end{array}\right)=b_1\cdots b_kb$.
\end{lemma}

\n\textit{Proof.}

\n (i) This follows since $\left(\begin{array}{cc}p&q\\r&s\end{array}\right) =\left(\begin{array}{c}\lambda\\\lambda'\end{array}\right) A (\gamma,\gamma')$.

\n (ii) Let $N=M(b_1,x,b_2)\cdots M(b_{k-1},x,b_k)$. Then $N=\left(\begin{array}{cc}b_1\cdots b_{k-1}&u\\0&b_2\cdots b_k\end{array}\right)$ and $(1,a)N=(b_1\cdots b_{k-1},u+ab_2\cdots b_k)$. Thus 
$$\lambda'=(1,a)N\left(\begin{array}{cc}b&0\\1&b_k\end{array}\right)=(bb_1\cdots b_{k-1}+u+ab_2\cdots b_k,ub_k+ab_2\cdots b_{k-1}b_k^2).$$ 
It follows that
$$det\left(\begin{array}{c}\lambda'\\ \lambda\end{array}\right)=(bb_1\cdots b_{k-1}+u+ab_2\cdots b_k)b_k-(ub_k+ab_2\cdots b_{k-1}b_k^2)=bb_1\cdots b_k.$$
\begin{flushright}
$\square$
\end{flushright}

\begin{lemma} \label{lemmapqrsb}
Let $A\in R^{2\times 2}$ and $a,b_1,...,b_k,b,c,c_1,...,c_l,d\in R$, $k,l\geq 1$. Define 
$$
\begin{array}{lll}
p&=&(1,b_k)A\left(\begin{array}{c}1\\ c_1\end{array}\right),\\
q&=&(1,b_k)A M(c,x,c_1)M(c_1,y,c_2)\cdots M(c_{l-1},y,c_l)  \left(\begin{array}{c}1\\ d\end{array}\right), \\
r&=&(1,a)M(b_1,x,b_2)\cdots M(b_{k-1},x,b_k)M(b_k,y,b) A \left(\begin{array}{c}1\\ c_1\end{array}\right), 
\\
s&=&(1,a)M(b_1,x,b_2)\cdots M(b_{k-1},x,b_k)M(b_k,y,b) AM(c,x,c_1)M(c_1,y,c_2)\cdots \\
&&M(c_{l-1},y,c_l)  \left(\begin{array}{c}1\\ d\end{array}\right).
\end{array}
$$
Then $det\left(\begin{array}{cc}p&q\\r&s\end{array}\right) =b_1\cdots b_k b c c_1\cdots c_l det(A)$.
\end{lemma}

\n\textit{Proof.}
Let $\lambda=(1,b_k)$, $\gamma=\left(\begin{array}{c}1\\ c_1\end{array}\right)$, $\gamma'=M(c,x,c_1)M(c_1,y,c_2)\cdots M(c_{l-1},y,c_l)  \left(\begin{array}{c}1\\ d\end{array}\right)$, $\lambda'=(1,a)M(b_1,x,b_2)\cdots M(b_{k-1},x,b_k)M(b_k,y,b)$. By Lemma \ref{lemmapqrsa}, $det\left(\begin{array}{cc}p&q\\r&s\end{array}\right) =det(A)det\left(\begin{array}{c}\lambda\\\lambda'\end{array}\right) det(\gamma,\gamma')$. By Lemma \ref{lemmapqrsa} again, $det\left(\begin{array}{c}\lambda\\\lambda'\end{array}\right)=-b_1\cdots b_kb$ and symmetrically, $det(\gamma,\gamma')=-cc_1\cdots c_l$, which ends the proof.
\begin{flushright}
$\square$
\end{flushright}

\b
\n\textbf{Proof of Theorem \ref{tiling-variables}.}
We prove that $t$ given by Eq. (\ref{tilingformula}) is an $SL_2$-tiling of the plane. It is enough to show that for any $(u,v)\in \mathbf Z^2$, the determinant of the matrix $\left ( \begin{array}{ll} t(u,v) & t(u,v+1)\\ t(u+1,v) & t(u+1,v+1) \end{array} \right )$ is equal to 1.

By inspection of the figure below, where $k,l\ge 1$ and $a,b_1,...,b_k,b,c,c_1,..,c_l,d$ are in $R$ and $w=x_1e_1\cdots e_{n-1}x_n$, $n\geq 0$, with $e_i\in R$ and $x_i\in\{x,y\}$,
$$
\begin{array}{cccccccccccc}
&&&&&&&&&&c_l&d \\
&&&&&&&&&&.     &|      \\
&&&&&&&&&&.     &|     \\
&&&&&&&&&&.      &|    \\
&&&&&&&&&&c_2 &|  \\
&&&&&&&&&c&c_1&|  \\
&&&&&&&&.&&|&| \\
&&&&&&&w&&&|&|  \\
&&&&&&.&&&&|&|     \\
&&&&&b&&&&&|&|  \\
b_1&.&.&.&b_{k-1}&b_k&-&-&-&-&(u,v)&(u,v+1) \\
a&-&-&-&-&-&-&-&-&-&(u+1,v)&(u+1,v+1)
\end{array}
$$
\n it is seen that the words associated to the four points $(u,v)$, $(u,v+1)$, $(u+1,v)$ and $(u+1,v+1)$ are  respectively of the forms $b_kybwcxc_1$, $b_kybwcxc_1yc_2\cdots yc_lxd$, $ayb_1x\cdots b_{k-1}xb_kybwcxc_1$ and $ayb_1x\cdots b_{k-1}xb_kybwcxc_1yc_2\cdots yc_lxd$.

Let $A=M(b,x_1,e_1)M(e_1,x_2,e_2)\cdots M(e_{n-1},x_n,c)$ and $D=be_1e_2\cdots e_{n-1}c$. Then define $p,q,r,s$ as in Lemma \ref{lemmapqrsb}. Thus we have $t(u,v)=\frac{p}{D}$, $t(u,v+1)=\frac{q}{Dc_1\cdots c_l}$, $t(u+1,v)=\frac{r}{b_1\cdots b_kD}$ and $t(+1,v+1)=\frac{s}{b_1\cdots b_kDc_1\cdots c_l}$. Therefore by Lemma \ref{lemmapqrsb}, 
$\Delta=t(u,v)t(u+1,v+1)-t(u,v+1)t(u+1,v)$
$=\frac{1}{b_1\cdots b_kD^2c_1\cdots c_l}(pq-rs)$ $=\frac{b_1\cdots b_kbcc_1\cdots c_l}{b_1\cdots b_kD^2c_1\cdots c_l} det(A)$. Now $det(A)=be_1^2e_2^2\cdots e_{n-1}^2c$ and $D^2=b^2e_1^2e_2^2\cdots e_{n-1}^2c^2$. Therefore $\Delta=1$.

\begin{flushright}
$\square$
\end{flushright}

\begin{corollary}
The sequences of a frise of type $\tilde A_m$, whose initial values are variables, are rational over the semiring of Laurent polynomials generated by these variables.
\end{corollary}

The proof is quite analogue to the proof in Section 7.1. Observe that if one follows the lines of the proof, one may recover the formula for the sequence associated to frise of the Kronecker quiver, as given in Section 2. In particular, one verifies that
$$M(a,x,b)M(b,y,a)=M,$$
with the notations in Section 2.

\subsection{Partial tilings and explicit formulas in case $A_n$}

We call {\em partial $SL_2$-tiling of the plane} a filling of a subset of $\mathbf Z^{2}$ such that each connected $2\times 2$ submatrix is of determinant 1. We construct some of these partial tilings. Note that our construction below is somewhat equivalent to the construction of frieze patterns of Coxeter \cite{Co} and \cite{Cacha}.

We begin with a construction which we call the {\em cross construction}: consider a word $w$, for example $w=aybycxdxexfxgyhyiyj$ (with $a,b,c,d,e,f,g,h,i,j$ either variables or set equal to 1) and its transpose $^tw=jxixhxgyfyeydycxbxa$. These two words determine two discrete paths in the plane; we represent each path by the sequence of variables of the word; we put the second one at the south-east of the first, separated by a cross, as shown below by asterisks, and then fill diagonally by 1's, italicized below.
$$
\begin{array}{cccccccccccc}
&&&&j&\it 1 \\
&&&&i&*&\it 1 \\
&&&&h&*&&\it 1 \\
c&d&e&f&g&*&&&\it 1 \\
b&&&&&*&&P&&\it 1 \\
a&&&&&*&&&&&\it 1 \\
\it 1&*&*&*&*&*&*&*&*&*&*&\it 1 \\
&\it 1&&&&*&&&&c&b&a \\
&&\it 1&&&*&&&&d \\
&&&\it 1&&*&&&&e \\
&&&&\it 1&*&&&&f \\
&&&&&\it 1&j&i&h&g 
\end{array}
$$

Then there exists a unique partial $SL_2$-tiling which fits into this figure. The proof of this result is a straightforward generalization of the proofs given so far. The only new thing is the definition of the word associated to each point in the region of the plane constructed above. Note that this region has naturally four components, separated by the cross (although they are not disjoint). 

For the north-west component, the word is defined by projection on the frontier, exactly as it has been done in Section 8.1. 

For the north-east component, one projects horizontally the point on the north-west frontier, and vertically on the south-east frontier: this give two words and the actual word is given by intersecting them. We do it on an example, which should be explicit enough. Consider the point denoted by $P$ above. Then the horizontal projection gives the point corresponding to the variable $b$ on the north-west frontier, and the vertical projection gives the point corrresponding to the variable $i$ on the south-east frontier. Then, the word associated to $P$ is by definition the word $bycxdxexfxgyhyi$; its denominator is, similarly to Section 8.1, the product of all the variables of the word, except the extreme ones. Thus, for $P$, it is $cdefgh$.

Then the value of the tiling at $P$ is 
$$
\frac{1}{cdefgh}(1,b)M(c,x,d)M(d,x,e)M(e,x,f)M(f,x,g)M(g,y,h) \left( \begin{array}{l}i\\1\end{array} \right),
$$
similarly to Eq. (\ref{tilingformulvariable}), except that the column matrix in the product is changed:

For the two other components, things are symmetric.

For the example, we show this tiling  below in the case where all variables are set equal to 1. 
$$
\begin{array}{cccccccccccc}
&&&&1&\it 1 \\
&&&&1&2&\it 1 \\
&&&&1&3&2&\it 1 \\
1&1&1&1&1&4&3&2&\it 1 \\
1&2&3&4&5&21&16&11&6&\it 1 \\
1&3&5&7&9&38&29&20&11&2&\it 1 \\
\it 1&4&7&10&13&55&42&29&16&3&2&\it 1 \\
&\it 1&2&3&4&17&13&9&5&1&1&1 \\
&&\it 1&2&3&13&10&7&4&1 \\
&&&\it 1&2&9&7&5&3&1 \\
&&&&\it 1&5&4&3&2&1 \\
&&&&&\it 1&1&1&1&1 
\end{array}
$$

The corresponding frise is then easily constructed; indeed, the word coding the south-east frontier is the transpose of the word coding the frontier; thus one continues the cross-construction with this new word. The tiling obtained has only to be repeated indefinitely; it is then seen to be periodic and its period is 2 + the length of the frontier $w$ (as follows already from the work of Fomin and Zelevinsky). If the frontier is an anti-palindrome, that is, $^tw=w$, then the period is the half of this number.

As a consequence of our results, we get the following well-known result (see
\cite{CK} for instance).

\begin{corollary} The Laurent phenomenon and the positivity conjecture hold for
cluster algebras of type $\tilde{A}_n$ or $A_m$ with acyclic initial clusters.
\end{corollary}

\section{Appendix: Cartan matrices of Dynkin and Euclidean types}
 
We represent Cartan matrices by their diagram. The Dynkin diagrams are in four infinite series, shown in Figure \ref{Dynkinseries}, where $m$ represents the number of vertices of each diagram, or one of the five exceptional ones, shown in Figure \ref{Dynkinexcept}.

\begin{figure}[htbp]
\begin{center}
\scalebox{.8}{\input{Dynkinseries.pstex_t}}
\caption{Infinite series of Dynkin diagrams}
\label{Dynkinseries}
\end{center}
\end{figure}

\begin{figure}[htbp]
\begin{center}
\scalebox{.8}{\input{Dynkinexcept.pstex_t}}
\caption{Exceptional Dynkin diagrams}
\label{Dynkinexcept}
\end{center}
\end{figure}

\begin{figure}[htbp]
\begin{center}
\scalebox{.8}{\input{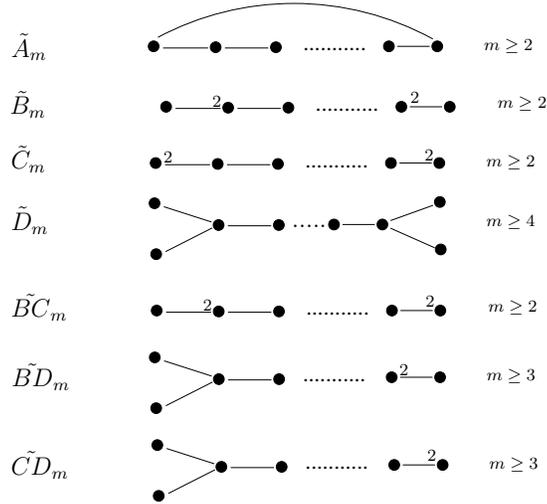}}
\caption{Infinite series of Euclidean diagrams}
\label{Euclideanseries}
\end{center}
\end{figure}

\begin{figure}[htbp]
\begin{center}
\scalebox{.8}{\input{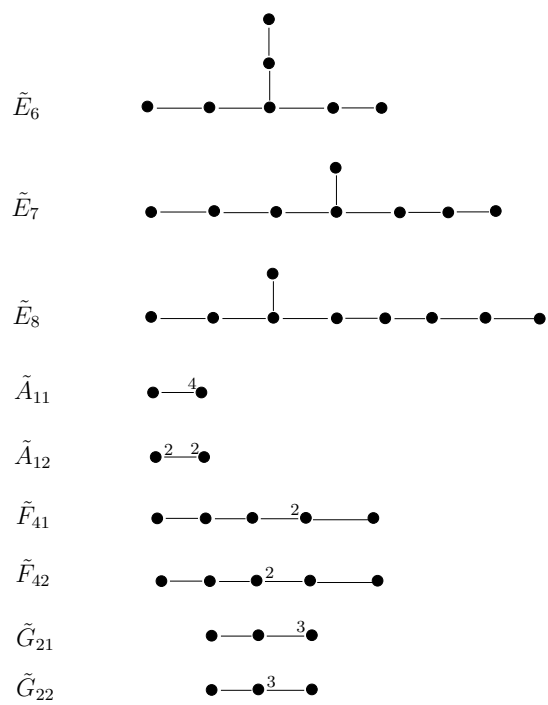}}
\caption{Exceptional Euclidean diagrams}
\label{Euclideanexcept}
\end{center}
\end{figure}

Furthermore, the Euclidean diagrams are in one of the seven infinites series shown in Figure \ref{Euclideanseries}, where $m$ is one less than the number of vertices of the diagram, or one of the nine exceptional ones, shown in Figure \ref{Euclideanexcept}.

\newcommand{\etalchar}[1]{$^{#1}$}

\end{document}